\documentclass{amsart}
\usepackage{amssymb}
\newtheorem{Theo}{Theorem}
\newtheorem{Lem}[Theo]{Lemma}
\newtheorem{Cor}[Theo]{Corollary}
\newcommand{\N}{\mathbb{N}}

\newcommand{\Q}{\mathbb{Q}}

\newcommand{\C}{\mathbb{C}}
\newcommand{\spec}{\mbox{spec}}
\begin{document}
\title{An additive property of almost periodic sets}
\author[J.-C. Schlage-Puchta]{Jan-Christoph Schlage-Puchta}
\begin{abstract}
We show that a set is almost periodic if and only if the associated exponential
sum is concentrated in the minor arcs. Hence binary additive problems involving
almost periodic sets can be solved using the circle method. This equivalence
is used to give simple proofs of theorems of J. Br\"udern.
\end{abstract}

\maketitle

A function $f:\N\rightarrow\C$ is called $\mathcal{B}^2$-almost periodic, if there
is a sequence of periodic functions $f_q$, such that
\[
\lim\limits_{q\rightarrow\infty}
\limsup\limits_{x\rightarrow\infty}\frac{1}{x}\sum_{n\leq x}|f(n)-f_q(n)|^2 = 0
\]
A set is called $\mathcal{B}^2$-almost periodic, if its characteristic function
is $\mathcal{B}^2$-almost periodic. The theory of almost periodicity is quite
rich, see e.g. \cite{SS}.

Let $\mathcal{N}$ be a set of integers. $\mathcal{N}$ is called distributed, if for any
$q$ and $a$, the density $f(q, a) = \lim\limits_{x\rightarrow\infty}
\frac{1}{x}\#\{n\leq x, n\in\mathcal{N}, n\equiv a\pmod{q}\}$ exists. A set
$\mathcal{N}$ is called extremal, if it has positive density $\rho$, is distributed,
and we have
\[
\frac{1}{\rho} = \sum_{q=1}^\infty\sum_{(a, q)=1}\left|\sum_{b=1}^q
\frac{f(q, a)}{\rho}e\left(\frac{ab}{q}\right)\right|^2
\]
Note that the sum over $a$ runs over residue classes prime to $q$, whereas the
sum over $b$ runs over all residue classes. Especially, if $q>1$ and $f(a, q)$
does not depend on $a$, the inner sum vanishes.
This definition is motivated by additive number theory. It turns out that
binary additive problems involving extremal sets can be solved using the circle
method. Define the major arcs $M(x, Q)=\bigcup\limits_{q\leq Q}\bigcup
\limits_{(a, q)=1}\left[\frac{a}{q}-\frac{Q}{x},\frac{a}{q}+\frac{Q}{x}\right]$
and the minor arcs $m(x, Q) = [0, 1]\setminus M(x, Q)$. Define $r(n)$ to be
the number of solutions of the equation $n=x+y$ with $x, y\in\mathcal{N}$. Then
we have
\[
r(n) = \int\limits_0^1 e(-n\theta) \left(\sum_{n\leq x\atop n\in\mathcal{N}}
e(\theta n)\right)^2 \;d\theta = \int\limits_{M(x, Q)} + \int\limits_{m(x, Q)}
\]
The integral over the major arcs can be evaluated whenever $\mathcal{N}$ is
distributed. Thus it remains to
bound the integral on the minor arcs. Hence one needs a nontrivial bound for
$S(\theta) = \sum\limits_{n\leq x\atop n\in\mathcal{N}}e(\theta n)$ on the minor
arcs. It turns out that this can be done for extremal sets. More precisely, we
have the following theorem. 

\begin{Theo}[Br\"udern]
Let $\mathcal{N}$ be a distributed set of positive density. Then the following is
equivalent.
\begin{enumerate}
\item $\mathcal{N}$ is extremal
\item As $Q$ and $x$ tend to infinity, we have $\int\limits_{m(x, Q)}
|S(\theta)|^2\;d\theta = o(x)$.
\end{enumerate}
\end{Theo}

Hence binary additive problems involving extremal sets can be solved. E.g. the
asymptotic number of representations of an integer as the sum of a $k$-free
and an $l$-free number can be computed\cite{FA}. Therefore a different
characterisation of extremal sets seems to be interesting. In this
note we prove the following theorem.

\begin{Theo}
The following two statements are equivalent:
\begin{enumerate}
\item $\mathcal{N}$ is extremal
\item $\mathcal{N}$ is $\mathcal{B}^2$-almost periodic
\end{enumerate}
\end{Theo}

Note that although additive questions involving almost periodic sets can be
deatl with in an elementary way, the theory of extremal sets gives better error
terms as shown in \cite{FA}.
From this one obtains the following corollaries.

\begin{Cor}
The intersection of extremal sets is extremal.
\end{Cor}

This was conjectured by J. Br\"udern\cite{B} and is the real motivation of the
present note. In the mean time J. Br\"udern gave a different proof (personal
communication).

\begin{Cor}
If $f(n) = \left\{\begin{array}{ll} 1 & \mbox{ if }
n\in\mathcal{N}\\0 & \mbox{ if } n\not\in\mathcal{N}\\\end{array}\right.$ is
multiplicative, and $\mathcal{N}$ has positive density, then $\mathcal{N}$ is
extremal.
\end{Cor}

This is theorem 1.4 in \cite{B}.

We will obtain theorem 2 as a corollary of a more general statement. To
formulate the next theorem, we have to introduce some notation.

We define $e_\beta(n):=e(\beta n)$.
Let $1\leq q<\infty$ be a real number. On the space of functions
$\N\rightarrow\C$ define a seminorm $\|f\|_q$ by
$\|f\|_q^q := \limsup\frac{1}{x}\sum_{n\leq x}|f(n)|^q$. Define the space
$\mathcal{B}$ of periodic functions, and the space $\mathcal{A}$ of trigonometric
polynomials $\sum_{\nu=1}^k a_\nu e_{\alpha_\nu}$ with $\alpha_\nu$ real.
Denote the closure of $\mathcal{B}$ with respect to $\|\cdot\|_q$ with $\mathcal{B}^q$,
and the closure of $\mathcal{A}$ with $\mathcal{A}^q$. For bounded
functions we have $f\in\mathcal{B}^q\Rightarrow f\in \mathcal{B}^{q'}$ for any
$q, q'$, and similar for $\mathcal{A}^q$. Define the scalar product
$\langle f, g\rangle$ by
\[
\langle f, g\rangle = \lim\limits_{x\rightarrow\infty}\frac{1}{x}\sum_{n\leq x}
f(n)\overline{g(n)}
\]
For $f, g\in\mathcal{A}^2$, this limit exists and defines a scalar product, which
induces the $\|\cdot\|_2$-seminorm. Hence we can apply the theory of Hilbert
spaces to obtain Fourier-series for almost periodic functions. Especially, if
$f\in\mathcal{A}^2$, for every $\beta$ the scalar product
$\langle f, e_\beta\rangle$ exists, and we define the Fouriercoefficient of
$f$ for $\beta$ to be this product. We define the spectrum of $f$ to be the
set of $\beta$, such that $\langle f, e_\beta\rangle\neq 0$ and write
$\spec(f)$ for this set. Now let $\alpha=(\alpha_n)$ be a sequence of numbers
from the interval $[0, 1]$. Define the major arcs $M_\alpha(x, Q)$ with
respect to this sequence by $M_\alpha(x, Q)=\bigcup\limits_{n\leq Q}
\left[\alpha_n-\frac{Q}{x}, \alpha_n+\frac{Q}{x}\right]$ and
$m_\alpha(x, Q)=[0, 1]\setminus M_\alpha(x, Q)$. A set $\mathcal{N}$ is
called $\alpha$-extremal if it has positive density, all Fouriercoefficients
of its characteristic function exist, and we have
$\int\limits_{m_\alpha(x, Q)} |S(\theta)|^2d\theta = o(x)$ for any $Q=Q(x)$,
tending to infinity with $x$.

Now we can state our main theorem.

\begin{Theo}\label{Haupt}
Let $\mathcal{N}$ be a set of integers with positive density, $f$ be the
characteristic function of $\mathcal{N}$. Let $\alpha = (\alpha_n)$ be some sequence
with $\alpha_n\in[0, 1)$. Then the following statements are equivalent:
\begin{enumerate}
\item $\mathcal{N}$ is $\alpha$-extremal
\item $f\in\mathcal{A}^2$, and $\spec(f)$ is contained in $\alpha$
\end{enumerate}
\end{Theo}

As theorem 1, this can be applied to additive questions.

\begin{Cor}\label{Wurzel}
Denote with $r(n)$ the number of positive integer solutions of the equation
$n=[\sqrt{2}a] + [\sqrt{3}b]$. Then we have for $n\rightarrow\infty$,
$r(n)\sim \frac{n}{\sqrt{6}}$.
\end{Cor}

Before we begin with proofs, we recall some facts about almost
periodic functions. All statements of this paragraph can be found in \cite{SS}.

\begin{Lem}
Let $f\in\mathcal{A}^2$ be a bounded. Then the  series
$\sum\limits_{\beta\in\mbox{\footnotesize spec}(f)}\langle
f, e_\beta\rangle e_\beta$ converges to $f$ with respect to the
$\|\cdot\|_2$-seminorm.
\end{Lem}

\begin{Lem}
The spectrum of any $f\in\mathcal{A}^2$ is countable. Furthermore,
$f\in\mathcal{B}^2$ if and only if $f\in\mathcal{A}^2$, and $\spec(f)\subseteq\Q$.
\end{Lem}

For proofs, see \cite{SS}, chapter VI.3.

\begin{Lem}
If $f, g\in\mathcal{A}^2$ are bounded, we get $fg\in\mathcal{A}^2$.
\end{Lem}

{\em Proof:} Since $fg$ is bounded, it suffices to show that $fg\in\mathcal{A}^1$.
This follows from \cite{SS}, theorem VI.2.3.

Further we will use the characterization of multiplicative almost periodic
functions.

\begin{Lem}
Let $f$ be a multiplicative function with mean value $M(f)\neq 0$ and
$|f(n)|\leq 1$ for all $n$. Then $f\in\mathcal{B}^2$ if and only if the following
series converge:
\begin{enumerate}
\item $S_1 = \sum_p \frac{1}{p}(f(p)-1)$
\item $S_2 = \sum_p \frac{1}{p}|f(p)-1|^2$
\end{enumerate}
\end{Lem}
{\em Proof:} This is a special case of \cite{SS}, Theorem VII.5.1.

First we show that theorem \ref{Haupt} implies theorem 2. Assume that theorem
\ref{Haupt} holds.
We claim that all the following statements are equivalent.
\begin{enumerate}
\item $\mathcal{N}$ is extremal
\item $\mathcal{N}$ is $\alpha$-extremal, where $\alpha$ is some denumeration of
the rational numbers in $[0, 1)$
\item $f\in\mathcal{A}^2$, and $\spec(f)$ is rational
\item $f\in\mathcal{B}^2$.
\end{enumerate}
The equivalence of 1. and 2. is obvious from the definition and theorem 1.
Note that the existence of the Fourier-coefficients is equivalent to the
fact that $\mathcal{N}$ is distributed. The
equivalence of 2. and 3. is given by theorem \ref{Haupt}, and the equivalence
of 3. and 4. is given by theorem 8. Hence 1. and 4. are equivalent, which
proves theorem 2.

Corollary 3 follows from theorem 2 and theorem 9, and corollary 4 follows from
theorem 2 and theorem 10, where the convergence of the series is implied by
the fact that $f$ takes values in $\{0, 1\}$, hence convergence of the series
is equivalent to the condition $M(f)\neq 0$.

To prove corollary \ref{Wurzel}, note that the condition
``$\exists a:[\sqrt{2}a] = k$''
is equivalent to the condition $[(k+1)/\sqrt{2}]-[k/\sqrt{2}]=1$. Replacing the
square brackets by approximating exponential polynomials, we see that the
characteristic function $f$ of the set $\{k|\exists a: [\sqrt{2}a] = k\}$ is
in $\mathcal{A}^2$ with $\spec{f} = k/\sqrt{2}\bmod{1}$, and the corresponding
statement is true for the set $\{k|\exists a: [\sqrt{3}a] = k\}$. Now
\[
r(n) = \int\limits_0^1 e(-n\theta)S_{\sqrt{2}}(\theta)S_{\sqrt{3}}(\theta)
d\theta
\]
where $S_{\sqrt{2}}(\theta)=\sum_{a\leq x/\sqrt{2}} e(\theta[\sqrt{2} a])$, and
$S_{\sqrt 3}$ is defined similar. Since we are not interested in an error term,
we choose $Q$ tending to infinity with $x$ sufficiently slowly.
Now if we define
$M_{\sqrt{2}}(x, Q):=\bigcup\limits_{|q|\leq Q}
\left[\left(\frac{q}{\sqrt{2}}\bmod{1}\right) - \frac{Q}{x},
\left(\frac{q}{\sqrt{2}}\bmod 1\right) + \frac{Q}{x}\right]$, and
$M_{\sqrt{3}}$ in the same way, we have for $x$
sufficiently large $M_{\sqrt{2}}(x, Q)\cap M_{\sqrt{3}}(x, Q)=
\left[\frac{-Q}{x},\frac{Q}{x}\right]$ since $\sqrt{2}$ and $\sqrt{3}$ are
linear independent over the rationals. This interval
contributes $\frac{n}{\sqrt{6}}+o(n)$ to the whole integral, hence it
suffices to estimate the remaining arcs. We have
\begin{eqnarray*}
\int\limits_{-\frac{Q}{x}}^{1-\frac{Q}{x}}
|S_{\sqrt{2}}(\theta)S_{\sqrt{3}}(\theta)\;d\theta & \leq &
\int\limits_{m_{\sqrt{2}}(x, Q)}|S_{\sqrt{2}}(\theta)S_{\sqrt{3}}(\theta)
d\theta\\
 && \quad + \int\limits_{m_{\sqrt{3}}(x, Q)}|S_{\sqrt{2}}(\theta)S_{\sqrt{3}}
(\theta) d\theta\\
 & \leq & \left(\int\limits_{m_{\sqrt{2}}(n, Q)}|S_{\sqrt{2}}(\theta)|^2\;
d\theta\cdot \int\limits_0^1|S_{\sqrt{3}}(\theta)|^2\;d\theta\right)^{1/2}\\
 &&\quad+ \left(\int\limits_0^1|S_{\sqrt{2}}(\theta)|^2\;
d\theta\cdot \int\limits_{m_{\sqrt{3}}(x, Q)}S_{\sqrt{3}}(\theta)|^2\;d\theta
\right)^{1/2}
\end{eqnarray*}
The second integral is $<n$, and the first integral is $o(n)$, by the
definition of $\alpha$-extremality. Hence the first summand is $o(n)$, and the
second summand can be dealt with similary.

Thus, it suffices to prove theorem \ref{Haupt}.

Assume that $f\in\mathcal{A}^2$, and let $(\alpha_n)$ be an enumeration of
$\spec(f)$. Choose $\epsilon>0$. Then there is some $N$,
such that $\|f-\sum_{\nu\leq Q}
\langle f, e_{\alpha_\nu}\rangle e_{\alpha_\nu}\|^2_2\leq \epsilon$.
Using the orthogonality of $e(\alpha n)$ and Parsevals equation we get
\begin{eqnarray*}
\int\limits_0^1\left|\sum_{n\leq x}\left(f(n) - \sum_{\nu\leq Q}
\langle f, e_{\alpha_\nu} \rangle e_{\alpha_\nu}\right)e(\theta n)
\right|^2\;d\theta
 & = & \|f-\sum_{\nu\leq N} \langle f,e_{\alpha_\nu}\rangle
e_{\alpha_\nu}\|^2_2 x\\
 & \leq & \epsilon x
\end{eqnarray*}
Hence to prove that $\int\limits _{m_\alpha(x, Q)}|S(\theta)|^2d\theta=o(x)$,
it suffices to prove this with $f$ replaced by some sufficiently long partial
sum of its Fourier-series. For if $g=\sum_{\nu\leq N}a_\nu e_{\alpha_\nu}$,
and $G$ is the corresponding exponential sum, we have
\[
\int\limits_{m_\alpha(x, Q)}|S(\theta)|^2d\theta \leq 2 \int\limits_{m_\alpha(x, Q)} |S(\theta)-G(\theta)|^2 d\theta + \int\limits_{m_\alpha(x, Q)}
|G(\theta)|^2 d\theta
\]
The first integral is $\leq\int\limits_0^1|S(\theta)-G(\theta)|^2d\theta =
\sum_{n\leq x}(f(n)-g(n))^2\leq \epsilon x$, thus it suffices to show that the
second integral is small, too. Now we have
\begin{eqnarray*}
\int\limits_{m_\alpha(x, Q)}\left|\left(\sum_{\nu\leq N}a_{\nu}
e_{\alpha_\nu}(n)\right) e(\theta n)\right|^2 d\theta & \leq &
\underbrace{\left(\sum_{\nu\leq N} |a_\nu|^2\right)}_{\leq 1}\cdot
\sum_{\nu\leq N}|a_{\nu}|\int\limits_{m_\alpha(x, Q)}\left|
\sum_{n\leq x}e((\theta-\alpha_\nu n)\right|^2 d\theta\\
 & \leq & 2N \int\limits_{Q/x}^{1-Q/x}\left|\sum_{n\leq x}
e(\theta n)\right|^2 d\theta\\
 & \leq & \frac{4Nx}{Q}
\end{eqnarray*}
Thus for any given $\epsilon>0$ we find some $N(\epsilon)$, such that 
\[
\int\limits_{m_\alpha(x, Q)} |S(\theta)|^2 d\theta \leq \epsilon x +
\frac{N(\epsilon)x}{Q}
\]
With $\epsilon\rightarrow 0$ and $Q\rightarrow\infty$, this becomes $o(x)$,
thus $\mathcal{N}$ is $\alpha$-extremal.

Now assume that $\mathcal{N}$ is an $\alpha$-extremal set of integers, and let
$\epsilon>0$. By the definition of $\alpha$-extremal, we get
$\int\limits_{m(x, Q)}|S(\theta)|^2d\theta = o(x)$, where
$M(x, Q)=\bigcup\limits_{q\leq Q}\bigcup\limits_{(a, q)=1}
\left[\frac{a}{q}-\frac{\omega(x)}{x}, \frac{a}{q}+\frac{\omega(x)}{x}\right]$,
$m(x, Q)=\left[\frac{\omega(x)}{x}, 1-\frac{\omega(x)}{x}\right]$, 
and $\omega(x)\nearrow\infty$ will be chosen later. Choose $Q$ such that for
all $x>x_0$ we have
$\int\limits_{m_\alpha(x, Q)} |S_x(\theta)|^2 \;d\theta < \epsilon x$.
Set $a_{\nu} = \lim\limits_{x\rightarrow\infty}\frac{1}{x}S(\alpha_\nu)$. Note
that this limit exists since all Fouriercoefficients of $f$ exist. Then
define $f_Q(n) = \sum_q\leq Q \sum_{\nu\leq Q} a_{\nu} e(-\alpha_\nu n)$.
Obviously,  $f_Q$ is a trigonometric polynomial, thus it
suffices to show that the sequence $f_Q$ approximates the characteristic
function of $\mathcal{N}$. Using orthogonality of $e(\alpha)$, we have
\begin{eqnarray*}
\sum_{n\leq x} |f(n) -f_Q(n)|^2 & = & \int\limits_0^1
|\sum_{n\leq x}(f(n) -f_Q(n))e(\theta n)|^2 d\theta\\
 & \leq & \int\limits_{M_\alpha(x, Q)} |\sum_{n\leq x}(f(n) -f_Q(n)) 
   e(\theta n)|^2 d\theta\\
 && + \int\limits_{m_\alpha(x, Q)} |\sum_{n\leq x}f(n)e(\theta n)|^2 d\theta + 
  \int\limits_{m_\alpha(x, Q)} |\sum_{n\leq x}f_Q(n)e(\theta n)|^2 d\theta\\
 & = & \int\nolimits_1 + \int\nolimits_2 +\int\nolimits_3
\end{eqnarray*}
The estimation of $\int\nolimits_2$ and $\int\nolimits_3$ is straightforward:
By assumption, $\int\nolimits_2\leq\epsilon x$, and the inequality
$\int\nolimits_3\ll\frac{x}{\omega(x)}$ can be deduced as above. Thus it
suffices to consider $\int\nolimits_1$. Since there are no more than $Q^2$
major arcs, it suffices to show that the integral taken over a single major arc
is $o(x)$. Hence we have to show that 
\[
\int\limits_{\alpha_\nu-\omega(x)/x}^{\alpha_\nu+\omega(x)/x} |\sum_{n\leq x}
(f(n)-f_Q(n)) e(\theta n)|^2 d\theta = o(x)
\]
Using partial summation the integral becomes
\[
\int\limits_{-\omega(x)/x}^{\omega(x)/x} \left|\sum_{n\leq x}(f(n)-f_Q(n))
e(\alpha_\nu n) -
\sum_{n\leq x}\sum_{k\leq n}(f(k)-f_Q(k))e(k\alpha_\nu)
(e(\theta(n+1))-e(\theta n))\right|^2\;d\theta
\]
Since $\sum_{n\leq x} f(n)e(\alpha_\nu n) \sim a_{\nu} x$ and
$\sum_{n\leq x} f_Q(n)e(\alpha_\nu n) = a_{\nu} x + O(1)$, there is some
function $\phi(x)$, tending monotonically to $\infty$, such that
$\sum_{n\leq x}(f(n)-f_Q(n))e(\alpha_\nu n)<\frac{x}{\phi(x)}$. Using this
function, the integral can be estimated by
\[
\int\limits_{\alpha_\nu-\omega(x)/x}^{\alpha_\nu+\omega(x)/x} \left|
\frac{x}{\phi(x)} +  \sum_{n\leq x} \frac{n}{\phi(n)}\theta\right|^2 d\theta <
\frac{2x \omega(x)^2} {\phi(x)^2} + \frac{2x\omega(x)^3}{\phi(\sqrt{x})^2} +
\omega^3(x)
\]
Putting the estimates for $\int\nolimits_i$ together, and summing over all the
$\ll Q^2$ major arcs, we get
\[
\sum_{n\leq x} |f(n) -f_Q(n)|^2 < \epsilon x + c\frac{x}{\omega(x)} + 
\frac{4xQ^2\omega(x)^3}{\phi(\sqrt{x})^2} + \omega^3(x)Q^2
\]
With $\omega(x) = \min(\phi(\sqrt{x})^{1/2}Q^{-2}, x^{1/4}Q^{-2})$, the right
hand side becomes $(\epsilon+o(1)) x$, since this is true for any
$\epsilon>0$, it is $o(x)$. Hence $f_Q\rightarrow f$ with respect to the
$\mathcal{A}^2$-norm, i.e. $f$ is $\mathcal{A}$-almost periodic, and $\spec{f}$ is
contained in $\alpha$.

\end{document}